\documentclass[12pt,draft]{amsart}

\usepackage{amssymb}
\usepackage{url}
\usepackage{amscd}
\usepackage{amsthm}
\usepackage{amsmath}

\newtheorem {Theorem}  {Theorem}

\def\:{\! :\!}

\begin{document}

\date{}

\title[A sum-division estimate of reals]
{A sum-division estimate of reals  }

\author{Liangpan Li and Jian Shen}

\address{Department of Mathematics, Shanghai Jiao Tong University,
Shanghai 200240, China $\&$ Department of Mathematics, Texas State
University, San Marcos, Texas 78666, USA }
 \email{liliangpan@yahoo.com.cn}

\address{Department of Mathematics, Texas State University, San Marcos,
Texas 78666, USA} \email{js48@txstate.edu}



\date{}

\begin{abstract}
Let $A$ be a finite set of positive real numbers. We present a
sum-division estimate:
\[|A+A|^2|A/A|\geq\frac{|A|^4}{4}.\]
\end{abstract}

\maketitle

\section{Introduction}

Let $A$ be a finite set of positive real numbers throughout. The
sum-set, product-set and ratio-set of $A$ are defined respectively
to be
\begin{align*}
A+A&=\{a+b:a,b\in A\}, \\
 AA&=\{ab:a,b\in A\},\\
A/A&=\{a/b:a,b\in A\}.
\end{align*}
A famous conjecture of Erd\"{o}s and Szemer\'{e}di \cite{Erdos}
asserts that for any $\alpha<2$, there exists a constant
$C_{\alpha}>0$ such that
\[\max\big\{|A+A|,|AA|\big\}\geq C_{\alpha}|A|^{\alpha}.\]
In a series of papers \cite{Chen,Elekes,Ford,Nathanson,S1,S2}, upper
bounds on $\alpha$ were found by many authors. One highlight in this
direction was a proof by Elekes \cite{Elekes} that $\alpha$ can be
taken $ \frac{5}{4}$. His argument utilized a clever application of
the Szemer\'{e}di-Trotter theorem on point-line incidences.
Recently, using the concept of multiplicative energy and an
ingenious geometric observation, Solymosi \cite{Solymosi08} obtained
if $A$ is not a singleton, then
\begin{equation}
\label{JointSolymosi}|A+A|^2|AA|\geq
\frac{|A|^{4}}{4\lceil\log_{2}|A|\rceil},
\end{equation}
which yields
\begin{equation}\label{MainSolymosi}
\max\big\{|A+A|,|AA|\big\}\geq
\frac{|A|^{4/3}}{2\lceil\log_{2}|A|\rceil^{1/3}}.
\end{equation}
One cannot completely drop  the  logarithmic term
 in (\ref{MainSolymosi}) since if we
choose $\widetilde{A}=\{1,2,\ldots,n\}$, then
\cite{Erdos1,Erdos2,Ford08,Te}
\begin{equation}\label{1955estimate}
|\widetilde{A}\widetilde{A}|=\frac{n^2}{(\ln n)^{\beta+o(1)}},\ \ \
\ \beta=1-\frac{1+\ln\ln2}{\ln 2}=0.0860713....
\end{equation}
There is a subtle difference between $|\widetilde{A}\widetilde{A}|$
 and $|\widetilde{A}/\widetilde{A}|$. In fact, Elekes and Ruzsa
 \cite{ElekesRuzsa} showed that  there exists a universal constant $\gamma>0$ such
 that
\begin{equation}|A+A|^6|A/A|\geq\gamma|A|^8,\end{equation}
which yields
 \[|\widetilde{A}/\widetilde{A}|\geq\frac{\gamma }{64}|\widetilde{A}|^2\]
 by choosing $A=\widetilde{A}$.
This leads to a natural question: how to give a joint estimate on
$|A+A|$ and $|A/A|$? It is not difficult to use the
Szemer\'{e}di-Trotter theorem on point-line incidences to show that
\begin{equation}|A+A||A/A|\geq C|A|^{5/2}\end{equation}
holds for some universal constant $C>0$. Besides, if we carefully
analyze Solymosi's proof of  (\ref{JointSolymosi}), then
\begin{equation}
\label{last}|A+A|^2|A/A|\geq
\frac{|A|^{4}}{4\lceil\log_{2}|A|\rceil}.
\end{equation}
The main purpose of this note is to drop  the term
$\lceil\log_{2}|A|\rceil$ in (\ref{last}):

\begin{Theorem}\label{ourresult}
Let $A$ be a finite set of positive real numbers. Then
\[ |A+A|^2 |A/A| \ge \frac {|A|^4}
4.\] This implies a sum-division estimate
\[\max\big\{|A+A|,|A/A|\big\}\geq \frac{|A|^{4/3}}{2}.\]
\end{Theorem}

There is an explanation on Theorem \ref{ourresult} in plane
geometry. View $\mathbb{R}^2$ naturally as the complex plane
$\mathbb{C}$. Given a finite set $A$ of positive real numbers,
denote by $\mbox{Rad}(A\times A)$ and $\mbox{Ang}(A\times A)$
respectively the radius-set and the angle-set    of $A\times A$.
Applying Theorem \ref{ourresult} with $\widehat{A}=\{a^2:a\in A\}$
yields
\[\max\big\{|\mbox{Rad}(A\times A)|,|\mbox{Ang}(A\times A)|\big\}\geq\frac{|A|^{4/3}}{2}.\]
This shows the angle-set  and the radius-set  of $A\times A$ cannot
be small simultaneously.

\section{Proof of the main result}

Suppose $|A/A|=y$ and $A/A=\big\{z_i\big\}_{i=1}^y$. Suppose $z_i$
has $m_i$ representations in $A\times A$, that is,
\[m_i=\Big|\big\{(a,b)\in A\times A:\frac{a}{b}=z_i\big\}\Big| \ \ \ \ (i=1,2,\ldots,y).\]
Without loss of generality we may order all $m_i$'s as follows:
\begin{equation}\label{order}m_1 \le m_2 \le \cdots \le m_{y}.\end{equation} Since $|A|^2 =
\sum_{i=1}^{y} m_i$, there exists a unique integer $k$, $1\leq k\leq
y$, such that
\[
\sum_{i=1}^{k-1} m_i < \frac {|A|^2 } 2 \le \sum_{i=1}^k m_i \le
km_k.
\]
 Hence
\begin{equation}\label{estimate1}
|A/A| \ge k \ge \frac {|A|^2} {2m_k}
\end{equation}
and \begin{equation}\label{estimate2} \sum _{i=k }^{y} m_i = \Big(
|A|^2 - \sum_{i=1}^{k-1} m_i \Big) \geq\frac { |A|^2}
2.\end{equation} By (\ref{order}) and  Solymosi's geometric
observation \cite{Solymosi08},
\begin{equation}\label{estimate3}
|A+A|^2=|(A\times A)+(A\times A)| \ge m_k \sum _{i=k }^{y} m_i.
\end{equation}
Multiplying (\ref{estimate1}), (\ref{estimate2}) and
(\ref{estimate3}) yields
\[ |A+A|^2 |A/A| \ge \frac {|A|^4}
4.\] This proves Theorem \ref{ourresult}.

\textbf{Remark}. Let $F_n=\{a/q:1\le a\le q\le n,(a,q)=1\}$ be the
set of Farey fractions of order $n$. It is well-known (\cite{Hua})
that $|F_n|\sim\frac{3}{\pi^2} n^2$ as $n\rightarrow\infty$.
Besides,
 it is not difficult to deduce from (\ref{1955estimate}) (see also  \cite{Ford08,Haynes}) that
 \[\max\{|F_n+F_n|,|F_n-F_n|,|F_nF_n|,|F_n/F_n|\}\leq C\frac{n^4}{(\ln n)^{\beta+o(1)}}\ \ \ \ (n\rightarrow\infty).\]
This shows generally one can not expect the estimate
\[\max\{|A+A|,|A/A|\}\asymp|A|^2\ \ \ \ (|A|\rightarrow\infty).\]
We thank Dimitris Koukoulopoulos for communicating this example
to us.

\section*{Acknowledgment}
Li's research was supported by the Mathematical Tianyuan Foundation
of China (No. 10826088). Shen's research was partially supported by
NSF (CNS 0835834) and Texas Higher Education Coordinating Board (ARP
003615-0039-2007).

\end{document}